\documentclass{article}
\usepackage[utf8]{inputenc}

\usepackage{amsmath, amssymb, latexsym}
\usepackage{amsfonts}
\usepackage[usenames]{color}

\newtheorem{theorem}{Theorem}

\newtheorem{lemma}[theorem]{Lemma}
\newtheorem{example}[theorem]{Example}
\newtheorem{remark}[theorem]{Remark}

\newcommand{\Z}{\mathbb{Z}}

\newcommand{\epr}{\hfill$\diamondsuit$\smallskip}
\newcommand{\pr}{\textit{Proof}. }

\newcommand{\aut}{\operatorname{Aut}}

\newcommand{\lnd}{\operatorname{LND}}

\title{On Isotropy group of Danielewski surfaces}
\author{
Rene Baltazar\\	
\texttt{renebaltazar.furg@gmail.com}\\
Marcelo Veloso\\
e-mail:
\texttt{veloso@ufsj.edu.br}}
\date{April 2, 2020}

\begin{document}

\maketitle

\begin{abstract}
In the present work we consider differential rings of the form $(\mathcal B,D)$ where $\mathcal B$ is a Danielewski surface and $D$ is a locally nilpotent derivation on $\mathcal B$. Influenced by several recent works, we describe the isotropy group  of a locally nilpotent derivation, $D$, on  Danielewski surfaces, in the cases $xy = \varphi(z)$, $x^ny=\varphi(Z)$, and $f(x)x = \varphi(z)$. 
\end{abstract}


\section{Introduction}
Throughout this paper $\mathbb{K}$  denotes  an algebraically closed field of characteristic zero.  An additive mapping $D: R \rightarrow R $ is said to be a \textbf{\textit{ derivation}}  of the $\mathbb{K}$-algebra $R$, if it satisfies the Leibniz rule: 
\[
D(ab)=aD(b)+D(a)b,
\]
for all $a,b \in R$. If $A$ is a subring of $R$, a derivation of $R$ satisfying $D(A)=0$ is called $D$ an $A$-derivation. The set of all derivations of $R$ is denoted by $Der(R)$ and by $Der_{A}(R)$,  the set of all $A$-derivations of $R$.

We denote by $Aut(R)$ the group of $\mathbb{K}$-automorphisms of $R$. Observe that  $Aut(R)$ acts on $Der_{\mathbb{K}}(R)$ by   conjugation
\[
(\rho, D) \rightarrow \rho^{-1}D\rho.
\]

Given a derivation $D\in Der_{\mathbb{K}}(R)$ its \textbf{\textit{isotropy subgroup}}, with respect to this group action, is defined by
\[
Aut(R)_D=\{\rho \in Aut(R) \mid \rho D= D\rho \}.
\]

In recent years, many authors have addressed properties of this specific subgroup. For instance, L. Mendes and I. Pan \cite{{Baltazar},{MendesPan}} proved that the subgroup o $\mathbb{K}$-automorphism of $\mathbb{K}[X, Y]$ which commutes with simple derivations is trivial. Moreover, the authors described the isotropy group of a Shamsuddin derivation and thus showed that $D$ is a simple Shamsuddin derivation if, and only if, $Aut(R)_D = \{id\}$. In the general case
of more variables, L. Bertoncello and D. Levcovitz \cite{BerLev}, obtained that the subgroup of $\mathbb{K}$-automorphisms of $\mathbb{K}[X_1 , \dots , X_n ]$ which commutes with simple
Shamsuddin derivations is trivial.

Surfaces defined by equations of the form $X^n Y = P (Z)$, have been
studied in many different contexts, see \cite{BiaVeldan, DaiDan, Dubo, lml3}. Of particular interest is their
connection to the cancellation problem. It is common to refer to these surfaces as Danielewski surfaces, since W. Danielewski \cite{Dani} was the first to observe their connection to the famous problem of cancellation.

Some authors (for instance \cite{BiaVeldan, DuboAdd}) consider a Danielewski surface (varieties) to be a surface defined, more generally, by equation
\[
f(X) Y = Q(X, Z) \mbox{ and } X_1^{m_1} X_2^{m_2} \cdots X_n^{m_n} Y = Q(X_1 , \dots , X_n , Z).
\]
In this paper, our goal is to study the isotropy group for Danielewski surfaces given equations of the form
\[
xy = \varphi(z), \quad x^n y = \varphi(z) \quad \text{and} \quad f (x)y = \varphi(z).
\] In the first part, we present some basics and well-known examples. The following sections are then dedicated to the study of the isotropy group in the above classes of Danielewski surfaces.

\section{Basic Concepts and Examples}

We denote the group of units of a ring $R$ by $R^*$. Furthermore we reserve capital letters $X,Y,Z$ to denote indeterminates over a ring $R$ and we denote the polynomial ring in $X,Y,Z$ over $R$ by $R[X,Y,Z]$. The polynomial ring in $n$-indeterminates over $R$ is also denoted by $R^{[n]}$. We reserve the symbol $\mathbb K$ to indicate a  field of characteristic zero and algebraically closed.

\begin{example}
Let $D=\partial_X$ be the partial derivative  with respect to $X$ of $\mathbb{K}[X,Y]$. Then
\[
Aut(\mathbb{K}[X,Y])_D=\{(X+p(Y), a+bY)  \mid a, b \in \mathbb{K} \mbox{ and }p(Y) \in \mathbb{K}[Y]\}.
\]	
See Example 25 of the $\cite{Baltese}$. \epr 
\end{example}

A derivation $D$ is said to be \emph{locally nilpotent} if for every $f\in R$ there is a $n\geq 0$ such that $D^n(f)=0$. By a result of Rentschler (\cite{Re1968}) we know that a locally nilpotent derivation in $\mathbb K[X,Y]$ is conjugate to a derivation of the form $u(X)\partial_Y$, where $u\in \mathbb K[X]$. Thus, the following example describes the isotropy group of a locally nilpotent derivation in two variables up to conjugation.

\begin{example}
	Let $D=f(X)\partial_Y$ be the derivation  of $\mathbb{K}[X,Y]$, where $f(X) \in \mathbb{K}[X]$. Then
	\[
	Aut(\mathbb{K}[X,Y])_D=\{(a+bX, p(X)+cY )  \mid a, b, c \in \mathbb{K}, \,bc\neq 0 \mbox{ and }p(X) \in \mathbb{K}[X]\}.
	\]	
\end{example}
	Let $\rho \in Aut(\mathbb{K}[X,Y])$, then 
	\[
    \rho(X)=a_0(X)+a_1(X)Y+\cdots +a_t(X)Y^t 
	\]
	and
	\[
	\rho(Y)=b_0(X)+b_1(X)Y+\cdots+ b_s(X)Y^s
	\]
	where $a_i(X), b_j(X) \in \mathbb{K}[X]$.
	
	\begin{enumerate}
		\item Since $\rho(D(X))=D(\rho(X))$ we have
		\[
		\begin{split}
		0=\rho(0)
		=&\rho(D(X))\\
		=&D(\rho(X))\\
		=&D(a_0(X)+a_1(X)Y+\cdots +a_t(X)Y^t)\\
		=&a_1(X)f(X)+2a_2(X)f(X)Y+\cdots+ta_t(X)f(X)Y^{t-1}
		\end{split}.
		\]
This implies that 
\[
a_i(X)=0 \mbox{ for } i=1,\ldots,t.
\]
Hence $\rho(X)=a_0(X)=q(X) \in \mathbb{K}[X]$. 
		
		\item As $\rho(D(Y))=D(\rho(Y))$ we have
		\[
		\begin{split}
		 f(\rho(X))
	    =&\rho(f(X))\\
		=&\rho(D(Y))\\
		=&D(\rho(Y))\\
		=&D(b_0(X)+b_1(X)Y+\cdots +b_s(X)Y^s)\\
		=&b_1(X)f(X)+2b_2(X)f(X)Y+\cdots+sb_s(X)f(X)Y^{s-1}
		\end{split}.
		\]
	This implies that 
	\[
	b_i(X)=0 \mbox{ for } i=2,\ldots,s.
	\]
	Hence $\rho(Y)= p(X) + r(X)Y \in \mathbb{K}[X]$. 
 \end{enumerate}
	Remember that since $\rho$ is an automorphism, then the determinant of the Jacobian matrix must be a  nonzero complex number. In other words,
	
	\[
	\begin{split}
	Det(Jac(\rho(X),\,\rho(Y)))
	=& Det\left(\begin{bmatrix}
	q'(X)& 0       \\
	p'(X)+r'(X)Y &r(X)  \\ 
	\end{bmatrix}\right)\\
	=& q'(X)r(X)\in  \mathbb{K}^*.
	\end{split}
	\]
	So $q'(X),\, r(X) \in \mathbb{K}^*$. Therefore, $\rho=(a+bX,\, p(X)+cY)$  with $a, b, c \in \mathbb{K}$, $bc\neq 0$ and $p(X) \in \mathbb{K}[X]$.

Recently, R. Baltazar and I. Pan \cite{BaltPan} have presented a criterion for deciding whether that the isotropy group of a derivation is an algebraic group of a polynomial differential ring in two variables. As a particular case, they showed that the isotropy group of a locally nilpotent derivation  never is an algebraic group.

The following observation was proposed by the referee of this work.

\begin{remark}
Given a locally nilpotent derivation $D$ of a domain $R$, the isotropy group $Aut(R)_D$ always contains as a subgroup an abelian group isomorphic to $({\rm ker(D)};+)$. Indeed, given any nonzero $f \in {\rm ker(D)}$, $fD$ is a locally nilpotent derivation whose associated 
exponential automorphism $\rho = e^{t fD}$ commutes with D. Since $R$ has transcendence degree 1 over $\rm ker(D)$, it follows that for every domain of dimension greater than or equal to $2$, $({\rm ker(D)};+)$ is an infinite dimensional unipotent group contained in $Aut(R)_D$. In particular, the fact that $Aut(R)_D$ is never an algebraic group for such domains is thus a standard general well-known fact, not specific
to case of polynomial rings considered in \cite{BaltPan}. Furthermore, in view of this observation, it seems that the real interesting
problem concerning isotropy groups of locally nilpotent $D$ derivations is to determine when $Aut(R)_D$ is strictly bigger that $({\rm ker(D)};+)$. Note that this correspondence is already partly mentioned in \cite{FiWa}. Here they define in particular the centralizer $\mathcal{C}_{pol}(\partial)$ of an locally nilpotent derivation $\partial$ of a polynomial ring
$\mathbb{C}[X]$ as the set of derivations $\delta$ which commutes with $\partial$ (in the sense that the
Lie bracket $[\partial$;$\delta]$ is equal to zero, Definition 2.2 p. 41). The same definition
makes sense for any algebra $R$ over a field, and the statement $2$. p. $42$ in their
paper asserting that $\mathcal{C}_{pol}(\partial)$ is a module over $\rm ker(\partial)$ implies through the correspondence $$ker(\partial) \ni f \mapsto e^{(f\partial)} \in Aut(R)_{\partial}$$ that $Aut(R)_{\partial}$ contains the additive group
 $({\rm ker(\partial)};+)$ as a subgroup.
\end{remark}
		
\section{The case $XY = \varphi(Z)$}

In this section we study  the isotropy  group of a class of Danielewski surfaces given by the most general equation 
\[
XY = \varphi(Z),
\]
is thus means that we study the ring 
\[
\mathcal B=\mathbb{K}[X,Y,Z]/(XY - \varphi(Z)),
\]
where  $\mathbb{K}$ is an algebraically closed field of characteristic zero, $X,Y$ and $Z$ are indeterminates over $\mathbb{K}$, $\varphi(Z)$ with $\deg(\varphi)>1$. Recall that we  can also  write $\mathcal B=\mathbb{K}[x,y,z]$, where $x,y$ and $z$ are the images of $X,Y$ and $Z$ under the canonical epimorphism $\mathbb{K}[X,Y,Z] \to \mathcal B$. 

 The next result describes a set of generators for the group of $\mathbb{K}$-automorphisms of the ring $\mathcal{B}$.

\begin{theorem}(\cite[Theorem]{lml1})
	\label{autxy}
	The group $\aut(\mathcal B)$ is generated by the following  $\mathbb{K}$-automorphisms: 
	
	\begin{enumerate}
		\item Hyperbolic rotations  $H(x)=\lambda x$, $H(y)= \lambda^{-1}y$ and $H(z)=z$, where $\lambda \, \in \mathbb{K}^*$.
		
		\item Involution $I(x)=y$, $I(y)=x$ and $I(z)=z$.
		
		\item Triangular $T(x)=x$, $T(y)=y +[\varphi(z+xh(x))-\varphi(z)]x^{-1}$, and $T(z)=z+xh(x)$, where $h \in \mathbb K^{[1]}$.
		 
		\item If $\varphi(z)=c(z-a)^d$ then rescalings $R(x)=x$, $R(y)=\lambda^dy$ and $R(z)=\lambda z +(\lambda -1)a$, if $\varphi(z)=z^d$; $\lambda \, \in \mathbb{K}^*$ should be added.
		
		\item If $\varphi(z)=(z-a)^iQ((z+a)^n)$ and $\mu \, \in \mathbb{K}$ is such that $\mu^n=1$, then a symmetry  $S(x)=x$, $S(y)=\mu^iy$ and $S(z)=\mu z+(\mu-1)a$ should be added.
	\end{enumerate}
\epr	
\end{theorem}

Let $G$ be the subgroup of $Aut(\mathcal{B})$ generated by the automorphisms  $I$ and  $T$, defined in the Theorem $\ref{autxy}$. The group $G$ is called the {\it tame subgroup} of $Aut(\mathcal B)$. In $\cite{{DaiDan}}$, with  the {\it Transitivity Theorem}, D. Daigle shows that  the action of  $G$ on $KLND(\mathcal B)$ is transitive.

The next result is a direct consequence of the Transitivity Theorem.

\begin{lemma}(\cite[Corollary 9.5.3]{Dailec})
Given any $D' \in LND(\mathcal B)$, there exists $\theta \in G$ such that $\theta \circ D' z\theta^{-1}=f(x)D$ for some $f(x) \in \mathbb K[x]$, where $D$ is the unique element of $LND(\mathcal B)$ satisfying 
\[
D(x)=0,\,\,\,D(y)=\varphi_z(z) \mbox{  and  } D(z)=x.
\]
\epr
\end{lemma}

Note that the previous lemma is a generalization of Rentschler's Theorem \cite{Re1968} (recall that $\mathbb K^2$ is a special case of Danielewski surface).

The following technical lemma is useful to prove the main result of this section, and the next section.

\begin{lemma} 
\label{gx}
Let  $g(X)$  be a nonzero polynomial of $\mathbb K[X]$. Suppose that
	\[
	g(X)=\lambda^ng(\lambda X), 
	\]
where $n\geq 1$. Then $\lambda$ is a root of unity. 
\end{lemma}
\pr
By assumption $g(X)$ is a nonzero polynomial. Then there exist $t\in \mathbb{N}$  and positive integers $0\leq n_0<\ldots<n_t$ such that 
\[
g(X)=\sum_{l=0}^{t} a_lX^{n_l},
\]
where $a_0,a_1,\ldots, a_{t} \in \mathbb K^*$. Since $g(X)=\lambda^ng(\lambda X)$ we have 

\[
\begin{split}
\sum_{l=0}^{t} a_lX^{n_l}=
 g(X)&=\lambda^ng(\lambda X)\\
     &=\lambda^n\sum_{l=0}^{t} a_l(\lambda X)^{n_l}\\
     &=\lambda^n\sum_{l=0}^{t} a_l\lambda^{n_l} X^{n_l}\\
     &=\sum_{l=0}^{t} a_l\lambda^{n+n_l} X^{n_l}.\\
\end{split}
\]
Thereby $a_{l}= \lambda^{n+n_l}a_{l}$ and thus $\lambda^{n+n_l}=1$ for $l=0,1,\ldots,t$. So $\lambda$, for each $l=0,1,\ldots,t$, is a $(n+n_l)$-th root of unity. Therefore, $\lambda$ is a root of unity. 
\epr

The main result this section is

\begin{theorem}
	Let $D$ be the nonzero derivation  of $\mathcal B$ defined by 
	\[
D(x)=0,\,\,\,D(y)=g(x)\varphi_z(z) \mbox{  and  } D(z)=g(x)x.
\] 
Then the isotropy group of $\mathcal{B}$, $Aut(\mathcal B)_D$, is generated by a finite cyclic group of form $$\{(\lambda x, \lambda^{-1}y, z)| \ \lambda \, \in \mathbb{K}^*\}$$  together with the triangular $\mathbb K$-automorphisms.
\end{theorem}
\pr
According to the notations of the previous theorem:
\begin{enumerate}
\item 
Note that
\[
H(D(z))=H(x)=g(\lambda x)\lambda x, \   D(H(z))=D(z)=g(x)x,
\]
and
\[
H(D(y))=g(\lambda x) \varphi_z(z), \ D(H(y))= \lambda^{-1}g(x)\varphi_z(z). 
\]
The previous provide that $g(x)=\lambda g(\lambda x)$. The equality $g(x)=\lambda g(\lambda x)$ implies that $\lambda$ is a root of unity, by Lemma $\ref{gx}$. Thus, the isotropy group contains the finite group generated by $\{(\lambda x, \lambda^{-1}y, z)\}$.

\item
Since that
\[
I(D(x))=0 \mbox{   and   } D(I(x))=D(y)=g(x)\varphi_z(z)
\]
we have that $ID \not= DI$. Then, $I \notin Aut(\mathcal B)_D$.

\item
Observe that
\[
T(D(y))=T(g(x)\varphi_z(z))=g(x)\varphi_z(z+xh(x))
\]
and
\[
\begin{split}
	D(T(y))
	=& D(y +[\varphi(z+xh(x))-\varphi(z)]x^{-1})\\
	=& D(y)+D(\varphi(z+xh(x))-\varphi(z))x^{-1}\\
	=& g(x)\varphi_z(z)+\varphi_{z}(z+xh(x))D(z)x^{-1}-\varphi_{z}(z)D(z)x^{-1}\\
	=& g(x)\varphi_z(z)+\varphi_{z}(z+xh(x))g(x)-g(x)\varphi_{z}(z)\\
	=& g(x)\varphi_z(z+xh(x)).
	\end{split}
\]
Thus we have $T(D(y))=D(T(y))$. In addition, it is easy to verify  that \linebreak $T(D(x))=D(T(x))$ and   $T(D(z))=D(T(z))$. Therefore,  $TD=DT$.

\item 
Since
\[
D(R(z))=D(\lambda z+(\lambda -1)a) = \lambda x g(x) \mbox{ and } R(D(z))= xg(x)
\]
we have $R(D(y))=D(R(y))$ if, only if, $\lambda=1$. 
\item 
Since
\[
D(S(z))=D(\mu z+(\mu -1)a)=\mu g(x)x \mbox{ and } S(D(z))=g(x)x
\]
we have $S(D(y))=D(S(y))$ if, only if, $\mu=1$.
\end{enumerate}
The result follows from items 1, 2, 3, 4 and 5.
\epr

\begin{remark}
The famous tame automorphisms are never in the isotropy group of a locally nilpotent derivation on $\mathcal B$.
\end{remark}

%{\bf Question}: The derivation  $F(y)=\varphi_z(z),\,\,\,E(x)=\varphi_z(z) \mbox{  and  } E(z)=x+y$ is a simple derivation?

\section{The case $X^nY = \varphi(Z)$}

In this section we find the isotropy group of a locally nilpotent derivation for a class of  Danielewski surfaces given by 
\[
\mathcal{B}=\mathbb{K}[X,Y,Z]/(X^nY - \varphi(Z)),
\]
where $n>1$ and $deg(\varphi)>1$.

There is a set of generators for the group $\aut(\mathcal B)$.

\begin{theorem}(\cite[Theorem 1]{lml3})
	\label{autxny}
	The group $\aut(\mathcal B)$ is generated by the following  $\mathbb K$-automorphisms.   
	
	\begin{enumerate}
		\item Hiperbolic rotations $H(x)=\lambda x$, $H(y)= \lambda^{-n}y$ and $H(z)=z$, where $\lambda \, \in K^*$.
		
		\item Triangular $T(x)=x$, $T(y)=y +[\varphi(z+x^nh(x))-\varphi(z)]x^{-n}$, and $T(z)=z+x^nh(x)$, where $h \in \mathbb K^{[1]}$.
		 
		\item If $\varphi(z)=z^d$ then the automorphisms  $R(x)=x$, $R(y)=\lambda^dy$ and $R(z)=\lambda z$, where $\lambda \, \in K^*$, should be added.
		
		\item If $\varphi(z)=z^i\phi(z^m)$ then the automorphisms $S(x)=x$, $S(y)=\mu^dy$ and $S(z)=\mu z$, where $\phi \in \mathbb K^{[1]}$   $\mu \, \in K$ and $\mu^m=1$, should be added.
	\end{enumerate}
	\epr
	\end{theorem}

The next lemma characterizes the locally nilpotent derivations of $\mathcal B$.

\begin{lemma}(\cite[Proposition]{lml3})
A derivation  $D$ of  $\mathcal B$ is locally nilpotent if and only if 
\[
D(x)=0,\,\,\,D(y)=g(x)\varphi_z(z) \mbox{  and  } D(z)=g(x)x^n,
\]
where $g(x)\in \mathbb{K}[x]$.
\epr
\end{lemma}

The following result  describes the isotropy group of a locally nilpotent derivation on $\mathcal{B}$.

\begin{theorem}
\label{isoxn}
Let $D$ be the derivation  of $\mathcal B$ defined by 
\[
D(x)=0,\,\,\,D(y)=g(x)\varphi_z(z) \mbox{  and  } D(z)=g(x)x^n.
\] 
Then the isotropy group of $\mathcal{B}$, $Aut(\mathcal B)_D$, is generated by the triangular automorphisms $T$ and the hyperbolic rotations $H$, where  $\lambda$ is a root of unity.
\end{theorem}
\pr
It is sufficient to verify which automorphisms, given by Theorem $\ref{autxny}$, commute with $D$.
\begin{enumerate}
    \item
It is clear that $H(D(x))=D(H(x))$. Now observe that 
\[
H(D(y))=g(\lambda x) \varphi_z(z) \mbox{ and } D(H(y))=\lambda^{-n}g(x)\varphi_z(z).
\]
Then
\[
g(\lambda x)=\lambda^{-n}g(x) \mbox{ or } g(x)=\lambda^{n}g(\lambda x).
\]
This implies, by Lemma $\ref{gx}$, that $\lambda$ is a root of unity. In the same way $H(D(z))=D(H(z))$ implies that  $\lambda$ is a root of unity. So $HD=DH$ implies that  $\lambda$ is a root of unity. 

\item 
It easy to check that $T(D(x))=D(T(x))  \mbox{ and }  T(D(z))=D(T(z))$.
Since 
$
T(D(y))=T(g(x)\varphi_z(z))=g(x)\varphi_z(T(z))=g(x)\varphi_z(z+x^nh(x))$,
and
\[
\begin{split}
	D(T(y))
	=& D(y +[\varphi(z+x^nh(x))-\varphi(z)]x^{-n})\\
	=& D(y)+D(\varphi(z+x^nh(x))x^{-n}-D(\varphi(z))x^{-n}\\
	=& g(x)\varphi_z(z) +\varphi_z(z+ x^nh(x))D(z+ x^nh(x))x^{-n}-\varphi_z(z)D(z)x^{-n} \\
	=& g(x)\varphi_z(z) +\varphi_z(z+ x^nh(x))g(x)x^nx^{-n}-\varphi_z(z)g(x)x^nx^{-n} \\
	=& g(x)\varphi_z(z) +\varphi_z(z+ x^nh(x))g(x)-\varphi_z(z)g(x) \\
	=& g(x)\varphi_z(z+ x^nh(x)).
	\end{split}
\]
we have $D(T(y))=g(x)\varphi_z(z+ x^nh(x))= T(D(y))$. Thus  $TD=DT$.
    
\item  If $\varphi(z)=z^d$ then the $\mathbb{K}$-automorphisms  $R(x)=x$, $R(y)=\lambda^dy$ and $R(z)=\lambda z$, where $\lambda \, \in \mathbb{K}^*$, should be added. It is a simple matter to show that $D(R(x))= R(D(x))$. Note that
\[
D(R(y))=D({\lambda}^d y)=\lambda^dD(y)=\lambda^dg(x)\varphi_z(z)=d\lambda^dg(x)z^{d-1},
\]
and
\[
R(D(y))=R(g(x)\varphi_z(z))=g(x)\varphi_z(R(z))= g(x)\varphi_z(\lambda z)=d\lambda^{d-1}g(x)z^{d-1}.
\]
Furthermore, 
\[
D(R(z))=D(\lambda z)=\lambda g(x)x^n \mbox{  and  } R(D(z))=R(g(x)x^n)=g(x)x^n.
\]
So $\lambda=1$ and finally we have $R=Id$.
    
\item The proceeding  to verify that $S=Id$ is analogous to the previous item. 
\end{enumerate}
\epr

Nowicki proved that there are simple derivations in two variables with an arbitrarily large Y-degree (see \cite{Now2008}) on $\mathbb{K} [X,Y]$. A natural question is to determine examples of derivations where the isotropy group is arbitrarily large. We note that the following derivation has this property on $\mathbb{K} [X,Y]$:

Let
$$d = X \partial_X + (Y^s + pX) \partial_Y \, , \qquad p \in \mathbb{K}^* \, , \quad s \geq 2$$
and write
$$\rho(X) = a_0(X) + a_1(X) Y + \ldots + a_r(X) Y^r \, , \quad a_r \neq 0$$
$$\rho(Y) = b_0(X) + b_1(X) Y + \ldots + b_t(X) Y^t \, , \quad b_t \neq 0$$
Then
$$\rho(d(X)) = \rho(X) =  a_0 + a_1 Y + \ldots + a_r Y^r$$
and
\begin{eqnarray*}
	d(\rho(X)) & = & d(a_0 + a_1 Y + \ldots + a_r Y^r) \\
	& = & X a_0'(X) + X a_1'(X) Y + \ldots + X a_r'(X) Y^r + \\
	& & \qquad \qquad \qquad + a_1(X) d(Y) + \ldots + r a_r(X) Y^{r-1} d(Y) \\
	& = & X a_0'(X) + X a_1'(X) Y + \ldots + X a_r'(X) Y^r + \\
	& & \qquad \qquad \qquad + a_1(X) (Y^s + pX) + \ldots + r a_r(X) Y^{r-1} (Y^s + pX)
\end{eqnarray*}

From $d(\rho(X)) = \rho(d(X))$ it follows that $a_r Y^r = r a_r(X) Y^{s + r-1}$ and so $r = 0$ (since $s \geq 2$). Thus
$$\rho(d(X)) = \rho(X) = a_0(X) = c_0 + c_1 X + \ldots + c_l X^l \, , \quad c_i \in  \, , \quad c_l \neq 0$$
and
\begin{eqnarray*}
	d(\rho(X)) & = & d(c_0 + c_1 X + \ldots + c_l X^l) \\
	& = & X c_1 + 2 c_2 X^2 + \ldots + l c_l X^l
\end{eqnarray*}
Again, from $d(\rho(X)) = \rho(d(X))$ it follows that $c_l X^l = l c_l X^l$, so $l = 1$, and $c_0 = 0$. Thus $\rho(X) = cX$ for some $c \in \mathbb{K}^*$.

\vspace{0.3cm}

Analogously,
\begin{eqnarray*}
	\rho(d(Y)) & = & \rho(Y^s + p X) = \rho(Y)^s + p \rho(X) \\
	& = & (b_0(X) + b_1(X) Y + \ldots + b_t(X) Y^t)^s + pcX
\end{eqnarray*}
and
\begin{eqnarray*}
	d(\rho(Y)) & = & d(b_0(X) + b_1(X) Y + \ldots + b_t(X) Y^t) \\
	& = & X b_0'(X) + X b_1'(X) Y + \ldots + X b_t'(X) Y^t + \\
	& & \qquad \qquad \qquad + b_1(X) d(Y) + \ldots + t b_t(X) Y^{t-1} d(Y) \\
	& = & X b_0'(X) + X b_1'(X) Y + \ldots + X b_t'(X) Y^t + \\
	& & \qquad \qquad \qquad + b_1(X) (Y^s + pX) + \ldots + t b_t(X) Y^{t-1} (Y^s + pX)
\end{eqnarray*}
From $d(\rho(Y)) = \rho(d(Y))$ it follows that $b_t(X) Y^{st} = t b_t(X)^s Y^{s+t-1}$, so $t = 1$ and $b_1(X)^{s-1} = 1$. In particular, $b_1 \in \mathbb{K}$ and $\rho(Y) = b_0(X) + b_1 X$. But then
\begin{eqnarray*}
	\rho(d(Y)) = d(\rho(Y)) & \implies & (b_0(X) + b_1 Y)^s + pcX = X b_0'(X) + b_1 (Y^s + pX) \\
	& \implies & s b_1^{s-1} Y^{s-1} b_0(X) = 0 \\
	& \implies & b_0(X) = 0.
\end{eqnarray*}
and
\begin{eqnarray*}
	\rho(d(Y)) = d(\rho(Y)) & \implies & b_1^s Y^s + pcX = b_1 Y^s + b_1 p X) \\
	& \implies & b_1 = c.
\end{eqnarray*}
Thus $\rho(X) = cX$ and $\rho(Y) = cY$, where $c \in \mathbb{K}^*$ with $c^{s-1} = 1$. Therefore
\[
Aut_d(\mathbb{K}[X,Y]) = \{(X,Y) \mapsto (c^i X,c^i Y) \mid c^{s-1} = 1, \, 0 \leq i \leq s-1 \} \cong \Z/(s-1).
\]

\begin{remark} Note that just as we offer an example of derivations where the isotropy group is arbitrarily large in $\mathbb{K}[X,Y]$, so there is a natural question about examples in $\mathcal B$. 
\end{remark}

\begin{example}
Let $\mathcal B$ be the ring $\mathbb K[x,y,z]$, where $x^ny = \varphi(z)$, $n\geq 2$. Let $D$ be a  derivation of $\mathcal B$, defined by
\[
D(x)=0, D(y)=\varphi_z(z) \mbox{ and } D(z)=x^n,
\]
and $\lambda$ be a primitive  $n$-th root of unity. Then, by Theorem $\ref{isoxn}$, $\aut(\mathcal B)_D$ is generated by Triangular automorphisms $T$ 
\[
 T(x)=x, \,  T(y)=y +[\varphi(z+x^2h(x))-\varphi(z)]x^{-2} \mbox{ and } T(z)=z+x^2h(x),
\] 
and  Hyperbolic rotations $H$,  where $h \in \mathbb K^{[1]}$ and  $\lambda^n=1$. Indeed, the automorphism
\[
H(x)=\lambda x, \,H(y)= \lambda^{-n}y \mbox{  and } H(z)=z
\]
satisfies $HD=DH$ if and only if $\lambda^n=1$, see Lemma $\ref{gx}$ with $g(x)=1$. 
\end{example}

\begin{example}
	Let $\mathcal B$ be the ring $\mathbb K[x,y,z]$, where $x^2y = \varphi(z)$. Let $D$ be a  derivation of $\mathcal B$, defined by
	\[
	D(x)=0, D(y)=(x+1)\varphi_z(z) \mbox{ and } D(z)=(x+1)x^2.
	\]
	Then the automorphism 
	\[
	H(x)=\lambda x, \,H(y)= \lambda^{-2}y \mbox{  and } H(z)=z
	\]
	satisfies $HD=DH$ if and only if $\lambda=1$. Therefore, by Theorem $\ref{isoxn}$,  $H=id$ and $\aut(\mathcal B)_D$ is generated by Triangular automorphisms $T$ such that 
	\[
	 T(x)=x, \,  T(y)=y +[\varphi(z+x^2h(x))-\varphi(z)]x^{-2} \mbox{ and } T(z)=z+x^2h(x),
	 \] 
	 where   $h \in \mathbb K^{[1]}$.
\end{example}

%{\bf Question}: Let $D$ and $E$ be derivations on $\mathcal B$ defined by

%\[
%D(x)=0,\,\,\,D(y)=\varphi_z(z) \mbox{  and  } D(z)=x^n;
%\]
%and
%\[
%E(y)=0,\,\,\,E(x)=\varphi_z(z) \mbox{  and  } E(z)=nyx^{n-1}.
%\] 
%There exists $\theta \in \aut(\mathcal B)$ such that $\theta^{-1}D\theta = E$?

\section{The case $f(X)Y = \varphi(Z)$}

In this section, we study  the isotropy  group  of any locally nilpotent derivation of the class of Danielewski surfaces, given by the  ring 
\[
\mathcal B=\mathbb K[X,Y,Z]/(f(X)Y - \varphi(Z)),
\]
where  $\mathbb K$ is an algebraically closed field of characteristic zero, $X,Y$ and $Z$ are indeterminates over $\mathbb K$, $\varphi(Z) \in \mathbb{K}[Z]$,  $m>1$, and $\deg(f)>1$. Recall that we  can also write $\mathcal B=\mathbb K[x,y,z]$, where $f(x)y = \varphi(z)$ and $x,y$ and $z$ are the images of $X,Y$ and $Z$ under the canonical epimorphism $\mathbb K[X,Y,Z] \to \mathcal B$.

Let $d$ be a derivation on the polynomial ring  $\mathbb K[X,Y,Z]$ such that $$d(X)=0, \quad d(Y)= \varphi_Z (Z), \quad \text{ and } \quad d(Z)=f(X).$$ It follows that $d$ is locally nilpotent, since it is easy to see that $d$ is triangular. Moreover, note that $d(f(X)Y - \varphi(Z))=0$ and, thus, $d$ induces a locally nilpotent derivation $\mathcal D$ on $\mathcal B$ given by 
$$\mathcal D(x)=0, \quad \mathcal D(y)=\varphi_{z}(z), \quad \text{ and } \quad \mathcal D(z)=f(x).$$ 

The next result characterizes the locally nilpotent nilpotents of the $\mathcal{B}$.

\begin{lemma}(\cite[Corolary 8]{BiaVeldan}) 
	\label{lndb} 
	If $\mathcal D$ is the derivation of $\mathcal B$ given by $\mathcal D(x)=0$, $\mathcal D(y)=\varphi_z(z)$ and $\mathcal D(z)=f(x)$, then $\lnd(\mathcal B)=\{h \mathcal D \mid h \in\mathbb K[x] \}$. \epr
\end{lemma}

There is a set of generators for the group of $\mathbb{K}$-automorphisms of the ring $\mathcal{B}$.

\begin{theorem}(\cite[Theorem 15]{BiaVeldan}) 
	\label{autb}
	The group $\aut(\mathcal B)$ is generated by the $\mathbb{K}$-automorphisms: 
	\begin{enumerate}
		\item $H$ defined by $H(x)=\lambda x$, where  $\lambda^s=1$, $H(y)=\lambda^jy$ and $H(z)=z$, if $f(X)=X^jh(X^s)$, with  $h \in \mathbb K^{[1]}$ having a nonzero root;
		\item $T$ defined by $T(x)=x$, $T(y)= y + [\varphi(z+h(x)f(x)) -\varphi(z)]f(x)^{-1}$, and  $T(z)=z + h(x)f(x)$,  where $h \in \mathbb K^{[1]}$;
		\item $R$ defined by $R(x)=x$, $R(y)=\lambda^dy$ and $R(z)=\lambda z$, if $\varphi(z)=z^d$;
		\item $S$ defined by  $S(x)=x$, $S(y)=\mu^iy$ and $S(z)=\mu z$ where $\mu^m=1$, if $\varphi(z)=z^i\phi(z^m)$, with $\phi\in \mathbb K^{[1]}$.
	\end{enumerate}
\epr
\end{theorem}

The following theorem is the main result of this section.

\begin{theorem}
	Let $D$ be a locally nilpotent derivation of $\mathcal B$. Then the isotropy group of $\mathcal{B}$, $	Aut(\mathcal B)_D$, is generated by a finite cyclic group of form $$\{(\lambda x, y, z)| \ \lambda \, \in \mathbb{K}^*\}$$  together with the triangular $\mathbb K$-automorphisms.
		\end{theorem}
\pr 
Let $D$ be a locally nilpotent derivation of $\mathcal B$. By Lemma $\ref{lndb}$ 
\[
D(x)=0,\,\, D(y)=g(x)\varphi_z(z) \mbox{ and } D(z)=g(x)f(x),
\]
for some $g(x) \in \mathbb K[x]$. The proof consists in determining which automorphisms, from Theorem  $\ref{autb}$, commute with $D$.

\begin{enumerate}
	\item
A trivial verification  shows that $TD=DT$ (similar to the case $X^nY-\varphi(Z)$).

\item 
Let $H$ be the automorphism  defined by $H(x)=\lambda x$, where  $\lambda^s=1$, \linebreak  $H(y)=\lambda^jy$ and $H(z)=z$, if $f(X)=X^jh(X^s)$, with  $h \in \mathbb K^{[1]}$ having a nonzero root. Note that
\[
H(D(x))=H(0)=0=D(\lambda x)=D(H(x)),
\]
\[
H(D(y))=g(\lambda x)\varphi_z(z) \mbox{ and } D(H(y))=\lambda^j\varphi_z(z)g(x)
\]
and 
\[
H(D(z))=g(\lambda x)f(\lambda x) \mbox{ and } D(H(z))=g(x)f(x).
\]
If $HD=DH$ we have $g(\lambda x)f(\lambda x)=g(x)f(x)$. It thus implies that $g(\lambda x)=g(x)$. And  since  that  $g(\lambda x)\varphi_z(z)=\lambda^j\varphi_z(z)g(x)$ we have $\varphi_z(z)=\lambda^j\varphi_z(z)$. So $\lambda^j=1$.
\item 
Let $R$ be defined by $R(x)=x$, $R(y)=\lambda^dy$ and $R(z)=\lambda z$, if $\varphi(z)=z^d$. Note that
\[
R(D(y))=g(x)\varphi_z(\lambda z) \mbox{ and } D(R(y))=\lambda^dg(x)\varphi_z(z).
\]
If $RD=DR$ we have $g(x)\varphi_z(\lambda z)=\lambda^dg(x)\varphi_z(z)$ and $\varphi_z(\lambda z)=\lambda^d\varphi_z(z)$. Since $\varphi_z(z)=dz^{d-1}$ we have $d{\lambda}^{d-1}z^{d-1}=d{\lambda}^{d}z^{d-1}$. This implies $\lambda=1$ and $R=id$. 

\item Let $S$ be defined by  $S(x)=x$, $S(y)=\mu^iy$ and $S(z)=\mu z$ where $\mu^m=1$, if $\varphi(z)=z^i\phi(z^m)$, with $\phi\in \mathbb K^{[1]}$. Note that
\[
S(D(z))=g(x)f(x) \mbox{ and } D(S(z))=\mu g(x)f(x).
\]
Then $g(x)f(x)=\mu g(x)f(x)$ which implies $\mu=1$. Therefore, if $DS=SD$ we have $S=Id$. 
\end{enumerate}
The result follows from items 1, 2, 3 and 4.
\epr

%\section{Questions}

%\begin{enumerate}

%\item What is the isotropy group of the triangular derivations?

%\item What is the isotropy group of a linear derivation of the Fermat ring?

%\item What is the isotropy group of a LND of the Danielewski surface?

%\end{enumerate}

\end{document}